\documentclass{article}
\usepackage{graphicx} 
\usepackage{Mystyle}

\title{\textbf{Applying Hurlbert's Linear Optimization Technique to Establish Bounds on Pebbling Numbers}}
\author{Lingwen Li}
\date{September 2024}

\begin{document}

\maketitle

\begin{abstract}\label{sec:abstract}

This paper explores the application of Hurlbert's Linear Optimization Technique to determine bounds on pebbling numbers. By applying Hurlbert's weight functions and optimization techniques, this paper derives upper bounds for specific graph families. The study offers a comprehensive analysis of these bounds and contributes to a broader understanding of pebbling numbers in graph theory. Additionally, I applied the weight function lemma to calculate upper bounds for certain graphs, such as the Petersen graph, the Bruhat graph, and trees.

\end{abstract}

\tableofcontents

\newpage

{\section{Introduction}}

A pebbling configuration is a set of pebbles distributed over the vertices of an undirected graph $G$. This refers to a function mapping $V(G)$ to $\mathbb{N}$. A move involves removing two pebbles from a vertex and placing one on an adjacent vertex. This is called a pebbling move. Our objective is to guarantee that following a succession of pebbling moves, a specific vertex will have at least one pebble. The minimum number of pebbles that ensures any vertex on the graph $G$ can be reached by one pebble is called the pebbling number.

Graph pebbling was first introduced by Lagarias and Saks and then the journey of pebbling started with Kleitman and Lemke. In 1989, Kleitman and Lemke presented a proof that strengthened a conjecture by Erdős and Lemke regarding zero-sum sequences in finite groups.

A zero-sum sequence is a sequence of elements from a finite group $G$ whose sum equals the identity element of $G$. The identity element is defined as the element in a group that, when the group operation is applied, leaves any other element unchanged. A simple example of identity is 0 in addition and 1 in multiplication. Typically, zero-sum sequence problems focus on determining whether a sequence can be rearranged or manipulated so that the sum of its elements equals zero. For example, consider an abelian group $G$ under addition with elements $\left\{ 4, 5, -9, 3, 1, 8, 10 \right\}$, a subsequence like $\left\{ 4, 5, -9\right\}$ sums to 0, which makes it a zero-sum subsequence.

Lemke and Kleitman \cite{LK} further refined the understanding of zero-sum sequences in finite groups by formulating more general conjectures, as illustrated by the following lemma:

\begin{Lemma}
Let $ \mathbb{Z}_n $ denote the cyclic group of order $ n $, and let $ |g| $ denote the order of an element $ g $ in the group to which it belongs. For any sequence $ (g_k)_{k=1,...,n} $ of $ n $ elements from $ \mathbb{Z}_n $, there exists a zero-sum subsequence $ (g_k)_{k \in \mathrm{K}} $ such that 
\[
\sum_{k \in \mathrm{K}} \frac{1}{|g_k|} \leq 1.
\]
\end{Lemma}

The concept of pebbling in graphs originated from Lagarias and Saks’s effort to provide an alternative proof to that of Kleitman and Lemke.  Simultaneously and independently, Chung \cite{Chung} proved the above lemma, originally stated in number-theoretic terms. Furthermore, Chung developed and implemented these ideas, offering a more natural and structural approach to delve deeper into zero-sum sequences by using pebbling graphs and introducing pebbling numbers for the first time \cite{EH}.

The pebbling graph model has its own significance in real-world network systems, such as transportation and supply chain models. Like these models, graph pebbling involves moving a resource, like a commodity, from one set of sources to a set of sinks (destinations of data flow) under specific constraints. For example, in network flow, it is necessary to optimize the flow of a commodity through a network, ensuring that the flow along edges is restricted and conserved through vertices. The ultimate aim is always to maximize the amount of commodity reaching the sinks. Similarly, the transportation model focuses on minimizing the total cost to meet supply and demand. However, the supply chain model seeks to satisfy demands with minimum inventory, often disregarding transportation costs. The rules of graph pebbling incorporate a unique constraint: the loss of the commodity itself during movement across an edge. This constraint mirrors real-world scenarios. Thus, graph pebbling models the challenge of meeting demands with minimal resources while accounting for inevitable losses, making it a valuable tool for us to research resource distribution in some complex systems. Since then, the focus of the study of graph pebbling has gradually shifted from the auxiliary proof of number theory to the exploration of its own properties.

In this paper, we will explore Glenn Hurlbert's perspective \cite{Hurlbert} on graph pebbling. While Fan Chung established the foundation for the study of pebbling numbers, Glenn Hurlbert significantly advanced the field with his contributions. Hurlbert's research, which includes papers on various aspects of pebbling such as optimal pebbling and cover pebbling, has been influential. Notably, one of his most well-known publications, "The Linear Optimization Technique," focused on effectively generating upper bounds for pebbling numbers. Based on this method, Hurlbert provided and proved the following lemma:

\begin{Lemma}

Let $C_n$ denote the cycle on $n$ vertices. Then for $k\geq 1$, we have:

\[\pi(C_{2k})=2^k\quad \text{and}\quad \pi(C_{2k+1})=2\left\lfloor \frac{2^{k+1}}{3} \right\rfloor+1,\] where $\pi$ represents the pebbling number.

\end{Lemma}

This paper is divided into the following sections: In Section \ref{Sec:Preliminaries}, we recall some fundamental definitions in graph theory as well as the definitions of some basic types of graphs, including paths, trees, cycles, and the Petersen graph, among others. In Section \ref{Sec:Pebbling Numbers}, we first introduce the rules of the Pebbling Game along with the definition of the pebbling number. Subsequently, we present some fundamental results and solve for the pebbling numbers, such as paths, complete graph, $C_5$, and the Petersen graph. Additionally, we list the pebbling number results for various graphs. In Section \ref{Sec:Linear Optimization Technique}, we introduce the weight function lemma and the technique of linear optimization, and we demonstrate its application in computing upper bounds for the pebbling numbers of various graphs. This includes direct calculations for upper bounds of the pebbling numbers of graphs such as the Petersen graph, Bruhat graph, and trees.

\section{Preliminaries}\label{Sec:Preliminaries}

We begin with some fundamental graph theory concepts. In this subsection, we are going to introduce some basic definitions and notations of graph theory. A solid understanding of these graph theory concepts will greatly aid in comprehending our research on pebbling.\newline

\begin{Definition}
A graph consists of two sets: $V$ and $E$. The set $V$ contains all the $vertices$, while the set $E$ contains all the $edges$. So, for a graph $G$, we have
\[V(G)=\left\{ v_{1},v_{2},v_{3}\cdots v_{n-1},v_{n}\right\},\]
and
\[E(G)=\left\{ \left\{ v_{i},v_{j}\right\}\right\} (i,j\in [ 1,n]).\]

For simplicity, we will denote an edge simply by $v_{i}v_{j}$ instead of $\left\{ v_{i},v_{j} \right\}$.\newline

 We define a graph $H$ to be a $subgraph$ of a graph $G$ if the vertex set of $H$ is equal to the vertex set of $G$, and the edge set of $H$ is equal to the edge set of $G$, as given by:
    \[V(H)\subseteq V(G)\text{ and }E(H)\subseteq E(G)\Leftrightarrow H\subseteq G.\]
    
\end{Definition}

\begin{Definition}
The $order$ of a graph is defined as the cardinality of its vertex set, and the $size$ is the cardinality of its edge set.\newline
\end{Definition}

\begin{Definition}
Given two vertices $u$ and $v$ in a graph G, we say $u$ and $v$ are $adjacent$ if and only if there is an edge connecting them. Otherwise, if $u$ and $v$ are $nonadjacent$, then they are not connected, as given by:

\[u,v\in V(G),\left\{u,v\right\}\in E(G)\Leftrightarrow u \text{ and } v \text{ are adjacent.} \]

We define the $neighbourhood$ of a vertex $v$, denoted by $N(v)$, to be the set of vertices adjacent to $v$, as given by:

\[N(v)=\left\{x\in V\mid vx\in E\right\}.\]

If an edge $e$ has a vertex $v$ as an end vertex, we say that $v$ is $incident$ with $e$.
\end{Definition}

\begin{Definition}

Based on the definition of vertices and edges, the $degree$ of a vertex $v$, denoted by $\deg(v)$, is defined as the number of edges incident to $v$.

For example, in the Petersen graph shown below, each vertex has a degree of 3, as illustrated:

\[\deg(v) = 3 \quad \text{for all } v \in V(P),\] where $V(P)$ represents the set of vertices in the Petersen graph (See Figure 1).

\begin{eqnarray*}
    \begin{tikzpicture}[scale=0.8]
    \tikzstyle{every node}=[thick,minimum size=2pt, inner sep=1pt]
    
    \node[circle, fill=black] (1) at (0,2){};
    \node[circle, fill=black] (2) at (-1.9,0.6){};
    \node[circle, fill=black] (3) at (-1.2,-1.6){};
    \node[circle, fill=black] (4) at (1.2,-1.6){};
    \node[circle, fill=black] (5) at (1.9,0.6){};

    \node[circle, fill=black] (6) at (0,0.8){};
    \node[circle, fill=black] (7) at (-0.6,0.5){};
    \node[circle, fill=black] (8) at (-0.6,-0.5){};
    \node[circle, fill=black] (9) at (0.6,-0.5){};
    \node[circle, fill=black] (10) at (0.6,0.5){};
    \draw (1) -- (2);
    \draw (2) -- (3);
    \draw (3) -- (4);
    \draw (4) -- (5);
    \draw (5) -- (1);
    \draw (6) -- (8);
    \draw (7) -- (9);
    \draw (8) -- (10);
    \draw (9) -- (6);
    \draw (10) -- (7);
    \draw (1) -- (6);
    \draw (2) -- (7);
    \draw (3) -- (8);
    \draw (4) -- (9);
    \draw (5) -- (10);
    
    \node[minimum size=0pt,inner sep=0pt,label=below: Figure 1: The Petersen Graph $P$.] (name) at (0,-2.5){};
    \end{tikzpicture}   
    \end{eqnarray*}
    
\end{Definition}

\begin{Definition}
Next, we will define some common elements in a graph.
\begin{itemize} 
    \item [(a)] A $walk$ in a graph is a sequence of (not necessarily distinct) vertices $v_{1},v_{2},...,v_{k}$ such that $v_{i}v_{i+1}\in E$ for $i=1,2,...,k-1$. Such a sequence is sometimes called a $v_{1}-v_{k}$ walk, and $v_{1}$ and $v_{k}$ are the end vertices of the walk.

    \item [(b)]A path is a walk in which all vertices are distinct.

    \item [(c)]A cycle, also called a closed path, is a path $v_1, v_2, \dots, v_k$ (where $k \geq 3$) that returns to the starting vertex, i.e., with the edge $v_k v_1$
\end{itemize}
\end{Definition}

\begin{exam}
    Using the fundamental definitions above, we can introduce some typical graphs:

\begin{itemize} 
    \item [(a)] A $complete$ $graph$ of order $n$, denoted by $K_{n}$($n\geq1$), is defined to be a graph in which every pair of distinct vertices is connected by a unique edge (See Figure 2).
   \begin{eqnarray*}
        \begin{tikzpicture}[scale=0.8]
         \tikzstyle{every node}=[thick,minimum size=4pt, inner sep=1pt]
         \node[circle, fill=black] (1) at (0,0.5){};
         \node[minimum size=0pt,inner sep=0pt,label=below:$K_1$] (name) at (0,0){};
        \end{tikzpicture}
          & \hspace{4mm}
             \begin{tikzpicture}[scale=0.5]
         \tikzstyle{every node}=[thick,minimum size=4pt, inner sep=1pt]
         \node[circle, fill=black] (1) at (-0.8,1){};
         \node[circle, fill=black] (2) at (0.8,1){};
          \draw (1)--(2);
         \node[minimum size=0pt,inner sep=0pt,label=below:$K_2$] (name) at (0,0){};
        \end{tikzpicture}
           \hspace{4mm}
    \begin{tikzpicture}[scale=0.8]
    \tikzstyle{every node}=[thick,minimum size=4pt, inner sep=1pt]
    \node[circle, fill=black] (1) at (-1,0){};
      \node[circle, fill=black] (2) at (1,0){};
      \node[circle, fill=black] (3) at (0,1.7){};
      \draw (1)--(2);
      \draw (3)--(1);
      \draw (3)--(2);
 \node[minimum size=0pt,inner sep=0pt,label=below:$K_3$] (name) at (0,-0.3){};
    \end{tikzpicture}
    \hspace{4mm}
        \begin{tikzpicture}[scale=0.8]
    \tikzstyle{every node}=[thick,minimum size=4pt, inner sep=1pt]
    \node[circle, fill=black] (1) at (-1,0){};
      \node[circle, fill=black] (2) at (1,0){};
      \node[circle, fill=black] (3) at (-1,2){};  
      \node[circle, fill=black] (4) at (1,2){};
      	\node(v-2) at(3.5,1)[minimum size=0pt, label=left:{$\cdots$}]{};
      \draw (1)--(2);
      \draw (3)--(1);
      \draw (3)--(2);
      \draw (4)--(1);
      \draw (4)--(2);
      \draw (4)--(3);
 \node[minimum size=0pt,inner sep=0pt,label=below:$K_4$] (name) at (0,-0.3){};
    \end{tikzpicture}\\
   & \text{Figure 2: Examples of Complete Graphs $K_n$ for $n = 1, 2, 3, 4$.}
    \end{eqnarray*}
    \item [(b)]A $cycle$ of order $n$, denoted by $C_{n}$, is a graph of a cycle on $n$ vertices (See Figure 3).
    \begin{eqnarray*}
    \begin{tikzpicture}[scale=0.8]
    \tikzstyle{every node}=[thick,minimum size=4pt, inner sep=1pt]
    \node[circle, fill=black] (1) at (-0.7,0){};
    \node[circle, fill=black] (2) at (0.5,0){};
    \node[circle, fill=black] (3) at (1,1.2){};
    \node[circle, fill=black] (4) at (0,2){};
    \node[circle, fill=black] (5) at (-1.1,1.2){};
    \draw (1)--(2);
    \draw (3)--(2);
    \draw (4)--(3);
    \draw (5)--(1);
    \draw[dotted,line width=1pt] (5)--(4);
 \node[minimum size=0pt,inner sep=0pt,label=below: Figure 3:  The Cycle Graph $C_n$ with $n$ Vertices.] (name) at (0,-0.8){};
    \end{tikzpicture}   
    \end{eqnarray*}
    \item [(c)]A $path$ of order $n$, denoted by $P_{n}$, is simply a graph of a path on $n$ vertices with $n-1$ length (See Figure 4).
    
    \begin{eqnarray*}
    \begin{tikzpicture}[scale=0.8]
    \tikzstyle{every node}=[thick,minimum size=4pt, inner sep=1pt]
    \node[circle, fill=black] (1) at (-1,0){};
      \node[circle, fill=black] (2) at (1,0){};
      \node[circle, fill=black] (3) at (3,0){};
      \node(4) at(4.7,0)[minimum size=0pt, label=left:{$\cdots$}]{};
       \node[circle, fill=black] (5) at (5,0){};
       \node[circle, fill=black] (6) at (7,0){};
      \draw (1)--(2);
      \draw (3)--(2);
      \draw (5)--(6);
 \node[minimum size=0pt,inner sep=0pt,label=below:Figure 4: The Path Graph $P_n$ with $n$ Vertices.] (name) at (3,-0.8){};
    \end{tikzpicture}
    \end{eqnarray*}
    \item [(d)]A $tree$ is a connected graph that contains no cycles (See Figure 5).  We define a $subtree$  of $T$ to be a subgraph  $T'$ of a tree $ T$. 
    \begin{eqnarray*}
    \begin{tikzpicture}[scale=0.8]
    \tikzstyle{every node}=[thick,minimum size=4pt, inner sep=1pt]
    \node[circle, fill=black] (1) at (0,2){};
    \node[circle, fill=black] (2) at (-1,1){};
    \node[circle, fill=black] (3) at (1,1){};
    \node[circle, fill=black] (4) at (-1.5,0){};
    \node[circle, fill=black] (5) at (-0.5,0){};
    \node[circle, fill=black] (6) at (0.5,0){};
    \node[circle, fill=black] (7) at (1.5,0){};
    \draw (1)--(2);
    \draw (1)--(3);
    \draw (2)--(4);
    \draw (2)--(5);
    \draw (3)--(6);
    \draw (3)--(7);
 \node[minimum size=0pt,inner sep=0pt,label=below: Figure 5: The Tree Graph T.] (name) at (0,-0.8){};
    \end{tikzpicture}   
    \end{eqnarray*}
    \item [(e)] \textit{The Petersen graph} (See Figure 6), named after Julius Petersen, is famous for its high level of symmetry and interesting properties.
    \begin{eqnarray*}
    \begin{tikzpicture}[scale=0.8]
    \tikzstyle{every node}=[thick,minimum size=2pt, inner sep=1pt]
    
    \node[circle, fill=black] (1) at (0,2){};
    \node[circle, fill=black] (2) at (-1.9,0.6){};
    \node[circle, fill=black] (3) at (-1.2,-1.6){};
    \node[circle, fill=black] (4) at (1.2,-1.6){};
    \node[circle, fill=black] (5) at (1.9,0.6){};

    \node[circle, fill=black] (6) at (0,0.8){};
    \node[circle, fill=black] (7) at (-0.6,0.5){};
    \node[circle, fill=black] (8) at (-0.6,-0.5){};
    \node[circle, fill=black] (9) at (0.6,-0.5){};
    \node[circle, fill=black] (10) at (0.6,0.5){};
    \draw (1) -- (2);
    \draw (2) -- (3);
    \draw (3) -- (4);
    \draw (4) -- (5);
    \draw (5) -- (1);
    \draw (6) -- (8);
    \draw (7) -- (9);
    \draw (8) -- (10);
    \draw (9) -- (6);
    \draw (10) -- (7);
    \draw (1) -- (6);
    \draw (2) -- (7);
    \draw (3) -- (8);
    \draw (4) -- (9);
    \draw (5) -- (10);
    
    \node[minimum size=0pt,inner sep=0pt,label=below: Figure 6: The Petersen Graph $P$.] (name) at (0,-2.5){};
    \end{tikzpicture}   
    \end{eqnarray*}

\end{itemize}
\end{exam}

\section{Pebbling Numbers}\label{Sec:Pebbling Numbers}

In this chapter, we begin by defining the pebbling number through a fun mini-game, making it easier and more enjoyable for readers to grasp the concept. After setting the stage with this playful introduction, we'll delve into some foundational theorems and lemmas, many of which can be easily proved using induction or contradiction.

\subsection{The Pebbling Game}

Referring to \href{https://www.people.vcu.edu/~ghurlbert/pebbling/pebb.html}{Hurlbert's Pebbling website page}, we can imagine a lively mini-game, called \textit{the Pebbling Game} on a graph $G$, between two friends, Peter and Connie. Peter is the one with the coins (pebbles) and a knack for strategy, while Connie is the mastermind behind placing pebbles (deciding the configuration) on the graph. Here’s how the game works:

\begin{itemize}
    \item \textbf{Peter the Pebble-Purchaser:} Peter buys $t$ pebbles with his hard-earned money and hands them over to Connie. Of course, Peter doesn't want to spend more than absolutely necessary. 
     \item \textbf{Connie the Configuration-Setter:} Connie takes the pebbles and cleverly distributes them across the vertices of the graph $G$. She also picks a special vertex $r$, called the root, that Peter needs to target.
    \item \textbf{Pebble Movement Rule:} Moving pebbles is not free here! To move a pebble from one vertex to its neighbour, Peter must use two pebbles: one as a toll to the graph and one that actually reaches the neighbour. Connie loves making Peter think hard about how to get pebbles where they need to go. To better understand the rule, we will use a simple graph to show the move (See Figure 7).

    \begin{eqnarray*}
    \begin{tikzpicture}[scale=0.5]
    \tikzstyle{every node}=[thick,minimum size=5pt, inner sep=1pt]
    
    \node[circle, fill=black] (1) at (-20,0){};
    \node[circle, fill=blue] (11) at (-20.8,0){};
    \node[circle, fill=blue] (12) at (-20.5,0.5){};
    \node[circle, fill=blue] (13) at (-20.5,-0.5){};
    \node[circle, fill=black] (2) at (-13,0){};
    \node[circle, fill=blue] (21) at (-12.5,0){};
    \node[circle, fill=black] (3) at (-9,0){};
    \node[circle, fill=blue] (31) at (-9.5,-0.5){};
    \node[circle, fill=red] (32) at (-5.75,0.5){};
    \node[circle, fill=red] (33) at (-5.25,0.5){};
    \node[circle, fill=black] (4) at (-2,0){};
    \node[circle, fill=blue] (41) at (-1.5,0){};
    \node[circle, fill=black] (5) at (2,0){};
    \node[circle, fill=blue] (51) at (1.5,-0.5){};
    \node[circle, fill=black] (6) at (9,0){};
    \node[circle, fill=blue] (62) at (9.5,0.5){};
    \node[circle, fill=blue] (63) at (9.5,-0.5){};

    \draw (1) -- (2);
    \draw[->, thick, green] (3) -- (4);
    \draw (5) -- (6);
    
    \node[minimum size=0pt,inner sep=0pt,label=below: Figure 7: Pebbling Rule with a Directed Move on a Path.] (name) at (-5.5,-2.5){};
    \end{tikzpicture}   
    \end{eqnarray*}

    \item \textbf{Winning the game:} Peter wins if he successfully places a pebble on the root vertex $r$ after a series of these tricky moves. If he can’t, Connie wins the Pebbling Game.

\end{itemize}
When Peter successfully places a pebble on $r$, we say the configuration $C$ is $r$-solvable. If no matter which vertex Connie picks as the root, Peter can always get a pebble there, the configuration $C$ is solvable. Now, if Peter, being a careful spender, buys just enough pebbles to guarantee he can win the game for a particular root $r$, we call this number $\pi(G,r)$. If he’s savvy enough to ensure victory no matter where Connie puts the root, the number of pebbles he needs is denoted by $\pi(G)$.\\

In the Petersen graph example (See Figure 8), blue dots near each vertex represent pebbles positioned at that vertex. In the left configuration, through a series of pebbling moves, it is guaranteed that a pebble can be placed on any vertex, making the configuration solvable. However, in the right configuration, it is impossible to move a pebble to any vertex in the inner circle, rendering the configuration unsolvable.

   \begin{eqnarray*}
    \begin{tikzpicture}[scale=0.7]
    \tikzstyle{every node}=[thick,minimum size=2pt, inner sep=1pt]
    
    \node[circle, fill=black] (1) at (-4,2){};
    \node[circle, fill=blue] (33) at (-3.8,2.2){};
    \node[circle, fill=black] (2) at (-5.9,0.6){};
    \node[circle, fill=black] (3) at (-5.2,-1.6){};
    \node[circle, fill=black] (4) at (-2.8,-1.6){};
    \node[circle, fill=black] (5) at (-2.1,0.6){};
    \node[circle, fill=blue] (30) at (-2,0.9){};
    \node[circle, fill=blue] (31) at (-1.8,0.6){};
    \node[circle, fill=blue] (32) at (-1.9,0.3){};

    \node[circle, fill=black] (6) at (-4,0.8){};
    \node[circle, fill=black] (7) at (-4.6,0.5){};
    \node[circle, fill=blue] (34) at (-4.4,0.7){};
    \node[circle, fill=blue] (35) at (-4.7,0.8){};
    \node[circle, fill=blue] (36) at (-4.8,0.3){};
    \node[circle, fill=black] (8) at (-4.6,-0.5){};
    \node[circle, fill=black] (9) at (-3.4,-0.5){};
    \node[circle, fill=black] (10) at (-3.4,0.5){};

    \node[circle, fill=black] (11) at (4,2){};
    \node[circle, fill=blue] (41) at (4.3,2.2){};
    \node[circle, fill=black] (12) at (2.1,0.6){};
    \node[circle, fill=blue] (42) at (1.9,0.8){};
    \node[circle, fill=black] (13) at (2.8,-1.6){};
    \node[circle, fill=blue] (43) at (2.6,-1.8){};
    \node[circle, fill=black] (14) at (5.2,-1.6){};
    \node[circle, fill=blue] (44) at (5.4,-1.8){};
    \node[circle, fill=black] (15) at (5.9,0.6){};
    \node[circle, fill=blue] (45) at (6.1,0.8){};

    \node[circle, fill=black] (16) at (4,0.8){};
    \node[circle, fill=black] (17) at (3.4,0.5){};
    \node[circle, fill=black] (18) at (3.4,-0.5){};
    \node[circle, fill=black] (19) at (4.6,-0.5){};
    \node[circle, fill=black] (20) at (4.6,0.5){};
    
    \draw (1) -- (2);
    \draw (2) -- (3);
    \draw (3) -- (4);
    \draw (4) -- (5);
    \draw (5) -- (1);
    \draw (6) -- (8);
    \draw (7) -- (9);
    \draw (8) -- (10);
    \draw (9) -- (6);
    \draw (10) -- (7);
    \draw (1) -- (6);
    \draw (2) -- (7);
    \draw (3) -- (8);
    \draw (4) -- (9);
    \draw (5) -- (10);

    \draw (11) -- (12);
    \draw (12) -- (13);
    \draw (13) -- (14);
    \draw (14) -- (15);
    \draw (15) -- (11);
    \draw (16) -- (18);
    \draw (17) -- (19);
    \draw (18) -- (20);
    \draw (19) -- (16);
    \draw (20) -- (17);
    \draw (11) -- (16);
    \draw (12) -- (17);
    \draw (13) -- (18);
    \draw (14) -- (19);
    \draw (15) -- (20);
    
    \node[minimum size=0pt,inner sep=0pt,label=below: Figure 8: Comparison of $r$-solvable (left) and $r$-unsolvable (right) configurations on the Petersen graph] (name) at (0,-2.5){};
    
    \end{tikzpicture}   
    \end{eqnarray*}

In short, $\pi(G)$ is Peter’s secret formula for winning the Pebbling Game every time, no matter where Connie tries to trip him up. For $\pi(G)$, we have:

\[\pi(G)=\text{max}_{r}\pi (G,r).\]

Alternatively, $\pi(G)$ is one more than the maximum number of pebbles that Peter needs for an unsolvable configuration with respect to any root vertex $r \in V(G)$.\\

Thus far, we have outlined the rules of the game. Moving forward, we will now discuss the formal mathematical definitions of the Pebbling Game.

We define a graph as follows: 

\[G=(V,E),\] where $V$ represents the set of vertices and $E$ represents the set of edges. Distributing a number of pebbles across the vertices yields a specific configuration $C$ on the graph, described by the function:

\[C:V(G)\to \mathbb{N}.\]

Here, $C(v)$ denotes the number of pebbles placed on vertex $r$. The total number of pebbles distributed across the entire graph is then given by:

\[\left| C \right|=\sum_{v\in V(G)}C(v).\]

Next, we'll perform some pebbling moves on $G$. For any edge $e=(u,v)\in E(G)$, a pebble can be moved from vertex $u$ to vertex $v$ if and only if there are at least two pebbles on $u$, i.e., $C(u)\ge 2$. During the move, only one pebble reaches its destination vertex $v$, while another pebble is paid as a toll and removed from the graph. Thus, we can obtain a new configuration after a $u-v$ move:

\[ C’(w)=\left\{
\begin{array}{rcl}
C(u)-2      &      & \text{if }w=u,\\
C(v)+1      &      & \text{if }w=v,\\
C(w)        &      & \text{otherwise}.
\end{array} \right. \]

Given a graph $G$ and an initial configuration $C$, we select a vertex to be the root $r$. Our objective is to transport one pebble to $r$ by repeatedly applying the moves defined earlier. The definitions of a solvable configuration and the function $\pi$ remain previously described. To reiterate, the pebbling number, denoted by $\pi(G)$, is the smallest natural number $n$ such that, for any vertex designated as the root and any initial configuration of $n$ pebbles on the graph, it is possible, through a sequence of pebbling moves, to place at least one pebble on the root vertex. Everything we discuss from now on will revolve around this pebbling number.

\subsection{Fundamental Theorems}

In the previous section, we introduced the basic concept of moving pebbles and defined the pebbling number. In this section, without using the Linear Optimization Technique, we will present and prove some related lemmas and examples to deepen our understanding of these concepts. We now turn to specific examples of graphs and their pebbling numbers. \\

\begin{Lemma}\label{Root Choice}

For a path $P_n$, Connie’s choice between $v_1$ or $v_n$ will make Peter need to buy the most pebbles to win the Pebbling Game.

\end{Lemma}

\begin{proof}

According to the definition of pebbling moves, the number of pebbles required increases with the distance a pebble must travel, since each move halves the number of pebbles that can be moved forward. Additionally, pebbling moves are only necessary in one direction along the path—once a pebble has passed a vertex, it does not need to be moved back. Therefore, the maximum number of pebbles is required if and only if the target vertex is either $v_1$ or $v_n$, as these endpoints represent the longest possible distance a pebble must traverse.
\end{proof}

\begin{Lemma}\label{Path PN}

For a path $P_{n}$, we have $\pi(P_{n})=2^{n-1}$.
    
\end{Lemma}

\begin{proof}

We prove this result by induction on the length $n$ of the path $P_n$. \\

\textbf{Base Case:} $n = 1$

For $n = 1$, the path $P_1$ consists of a single vertex. Clearly, the pebbling number $\pi(P_1) = 1$, which satisfies the formula $2^{1-1} = 1$. Therefore, the base case $n = 1$ holds.\\

\textbf{Base Case:} $n = 2$

Now, consider the case when $n = 2$. The path $P_2$ consists of two vertices, which we denote as $v_1$ and $v_2$.

\begin{itemize}
    \item \textbf{Case 1:} If we place exactly one pebble on each vertex, no moves are required to ensure that each vertex has at least one pebble. Hence, this configuration satisfies the condition.
    
    \item \textbf{Case 2:} If both pebbles are placed on a single vertex, say $v_1$, and the target is $v_2$, then we can move one of the pebbles from $v_1$ to $v_2$, thus ensuring that $v_2$ has at least one pebble.

    \item \textbf{Infeasibility:} Here we also need to prove one pebble is not adequate. If only one pebble is placed on $V_1$, it is impossible to move it to $V_2$, meaning we cannot guarantee a configuration with at least one pebble on each vertex.
\end{itemize}

\textbf{Induction Step:}

Assume that for a path $P_k$, the pebbling number is $\pi(P_k) = 2^{k-1}$. We need to show that $\pi(P_{k+1}) = 2^k$. Consider the path $P_{k+1}$. According to Lemma \ref{Root Choice}, selecting either $v_{1}$ or $v_{k+1}$ necessitates having the maximum number of pebbles. Therefore, we will assume that we aim to move a pebble to the first vertex $v_1$. To do this, the vertex $v_2$ must contain at least 2 pebbles. By the inductive hypothesis, the remaining path from $v_2$ to $v_{k+1}$ (which is isomorphic to $P_k$) requires at most $2^{k-1}$ pebbles somewhere on $P_{k+1}$ to ensure that one pebble reaches $v_2$. Therefore, to guarantee a pebble reaches $v_1$, we need $2 \times 2^{k-1} = 2^k$ pebbles on $v_{k+1}$. Hence, $\pi(P_{k+1}) = 2^k$, completing the induction. 
\end{proof}

Moreover, we can expand our exploration to include pebbling numbers for various families of graphs, which often exhibit different characteristics and complexities. For instance, we have seen that the pebbling number for a path $P_n$ is given by $2^{n-1}$. To illustrate this further, we will consider some graph families such as the complete graph $K_n$, where the pebbling number is known to be exactly $n$. 

It's easy to guess that Peter can't buy too few pebbles if he wants to win the Pebbling Game.

\begin{Lemma} \label{pi bound}
For a graph $G$ with $n$ vertices, we have:

\[\pi(G)\geq n.\]

\end{Lemma}

\begin{proof}

To prove that the pebbling number of $G$ is greater than $n$, the number of vertices in $G$, we can demonstrate that with only $n-1$ pebbles, there exists a configuration that is $r$-unsolvable for some root vertex $r \in V(G)$. Specifically, we can distribute one pebble to each of the $n-1$ vertices, leaving the root vertex without any pebbles. Since it is impossible to move a pebble to the root vertex from this configuration, this proves that the pebbling number of $G$ must be equal or larger than its number of vertices.

\end{proof}

Due to the connectivity properties of complete graphs, Peter may not need “too many” pebbles to win the pebbling game on complete graphs.
\begin{exam}
     $\pi(K_n)=n$, where $K_n$ denotes the complete graph on $n$ vertices.
     
\end{exam}
\begin{proof}
    Let $K_n$ be a complete graph with  vertices $v_1,\ldots,v_n$. Let $C:V(G)\to \mathbb{N}$ be a  configuration on the graph $K_n$, for arbitrary $k =1,2,\ldots,n$. If $C(v_k)\neq0$,  $C$ is solvable.  If $C(v_k)=0$, by the Pigeonhole Principle, there exists a vertex $v_j$ with more than one pebble. Thus, a new configuration after $v_j-v_k$ move is solvable. Therefore, we prove that $\pi(K_n,v_k)\leq n$, for all $k =1,2,\ldots,n$. Thus, $\pi(K_n)\leq n$. By lemma~\ref{pi bound}, we have $\pi(K_n)=n$.

\end{proof}
   
\begin{exam}\cite{Hurlbert99}\label{PF}
      $\pi(P)=10$, where $P$ denotes the Petersen graph.
\end{exam}

\begin{proof}

Recall that the Petersen graph consists of two cycles, each with 5 vertices (See Figure 9).

\begin{eqnarray*}
    \begin{tikzpicture}[scale=0.7]
    \tikzstyle{every node}=[thick,minimum size=2pt, inner sep=1pt]
    
    \node[circle, fill=red] (1) at (0,2){};
    \node[circle, fill=red] (2) at (-1.9,0.6){};
    \node[circle, fill=red] (3) at (-1.2,-1.6){};
    \node[circle, fill=red] (4) at (1.2,-1.6){};
    \node[circle, fill=red] (5) at (1.9,0.6){};

    \node[circle, fill=blue] (6) at (0,0.8){};
    \node[circle, fill=blue] (7) at (-0.6,0.5){};
    \node[circle, fill=blue] (8) at (-0.6,-0.5){};
    \node[circle, fill=blue] (9) at (0.6,-0.5){};
    \node[circle, fill=blue] (10) at (0.6,0.5){};
    
    \draw[red] (1) -- (2);
    \draw[red] (2) -- (3);
    \draw[red] (3) -- (4);
    \draw[red] (4) -- (5);
    \draw[red] (5) -- (1);
    
    \draw[blue] (6) -- (8);
    \draw[blue] (7) -- (9);
    \draw[blue] (8) -- (10);
    \draw[blue] (9) -- (6);
    \draw[blue] (10) -- (7);
    
    \draw (1) -- (6);
    \draw (2) -- (7);
    \draw (3) -- (8);
    \draw (4) -- (9);
    \draw (5) -- (10);
    
    \node[minimum size=0pt,inner sep=0pt,label=below: Figure 9: The Petersen Graph Highlighting Two Cyclic Structures.] (name) at (0,-2.5){};
    \end{tikzpicture}
\end{eqnarray*}

Thus, rather than finding the pebbling number of the Petersen graph directly, we will first consider the pebbling number of a 5-vertex cycle, $C_5$.

Accordingly, based on Lemma \ref{pi bound}, we obtain:

\[\pi(C_5)\geq 5.\]

\begin{Lemma}
For a cycle $C_5$, its pebbling number satisfies:

\[\pi(C_5) \leq 5.\]
\end{Lemma}

\begin{proof}

For a cycle $C_5$ and its root $r$, we can always cut the cycle into two separate paths both including the root $r$ (see Figure 10).

\begin{eqnarray*}
    \begin{tikzpicture}[scale=0.7]
    \tikzstyle{every node}=[thick,minimum size=2pt, inner sep=1pt]
    
    \node[circle, fill=red] (1) at (0,2){};
    \node[circle, fill=black] (2) at (-1.9,0.6){};
    \node[circle, fill=black] (3) at (-1.2,-1.6){};
    \node[circle, fill=black] (4) at (1.2,-1.6){};
    \node[circle, fill=black] (5) at (1.9,0.6){};
    
    \draw[thick][blue] (1) -- (2);
    \draw[thick][blue] (2) -- (3);
    \draw[thick][black] (3) -- (4);
    \draw[thick][green] (4) -- (5);
    \draw[thick][green] (5) -- (1);

    \node[minimum size=0pt,inner sep=0pt,label=below: Figure 10: The Cycle Graph $C_5$ with a Red Root and Two Path Graphs $P_3$.] (name) at (0,-2.5){};
    \end{tikzpicture}
\end{eqnarray*}

According to Lemma \ref{Path PN}, we have the pebbling number formula for a path, as given by:

\[\pi(P_n)=2^{n-1}.\]

Therefore, for the path $P_3$, the pebbling number is $2^2 = 4$. Given that we have 5 pebbles, they can be distributed across the two paths in three possible ways: $0+5$, $1+4$, and $2+3$.

In the cases of the $0+5$ and $1+4$ distributions, it’s straightforward to verify that the configuration is $r$-solvable. For the $2+3$ configuration, if even one pebble is placed on the middle vertex of $P_3$, the configuration remains $r$-solvable. Thus, the only potentially challenging scenario is when all three pebbles are placed on the end vertex (which is not the root). Similarly, in this configuration, the remaining 2 pebbles on the other $P_3$ can also be distributed in different ways. If both pebbles are on the middle vertex, the configuration is $r$-solvable. If one pebble is on the middle vertex and the other on the end vertex (which is not the root), or if both pebbles are on the end vertex (which is not the root), we need to examine the entire graph. In each of these cases, it can be shown that the configuration remains $r$-solvable (see Figures 11 and 12).

\begin{eqnarray*}
    \begin{tikzpicture}[scale=0.5]
    \tikzstyle{every node}=[circle, draw, thick, minimum size=10pt, inner sep=0pt]
    
    \node[fill=red] (1) at (-10,0){};
    \node[fill=white] (2) at (-8,0){0};
    \node[fill=white] (3) at (-6,0){3};
    \node[fill=white] (4) at (-4,0){2};
    \node[fill=white] (5) at (-2,0){0};

    \node[fill=red] (6) at (2,0){};
    \node[fill=white] (7) at (4,0){0};
    \node[fill=white] (8) at (6,0){4};
    \node[fill=white] (9) at (8,0){0};
    \node[fill=white] (10) at (10,0){0};

    \draw[green, thick] (1) -- (2);
    \draw[green, thick] (2) -- (3); 
    \draw[black, thick] (3) -- (4);
    \draw[blue, thick] (4) -- (5);

    \draw[green, thick] (6) -- (7);
    \draw[green, thick] (7) -- (8); 
    \draw[black, thick] (8) -- (9);
    \draw[blue, thick] (9) -- (10);

    \draw[->, blue, thick, bend left=45] (5) to (1);
    \draw[->, blue, thick, bend left=45] (10) to (6);

    \draw[->, thick] (-1,0) -- (1,0);
    \node[minimum size=0pt,inner sep=0pt,label=below: Figure 11] (name) at (0,-0.7){};
    \end{tikzpicture}
\end{eqnarray*}

\begin{eqnarray*}
    \begin{tikzpicture}[scale=0.5]
    \tikzstyle{every node}=[circle, draw, thick, minimum size=10pt, inner sep=0pt]
    
    \node[fill=red] (1) at (-10,0){};
    \node[fill=white] (2) at (-8,0){0};
    \node[fill=white] (3) at (-6,0){3};
    \node[fill=white] (4) at (-4,0){1};
    \node[fill=white] (5) at (-2,0){1};

    \node[fill=red] (6) at (2,0){};
    \node[fill=white] (7) at (4,0){0};
    \node[fill=white] (8) at (6,0){1};
    \node[fill=white] (9) at (8,0){2};
    \node[fill=white] (10) at (10,0){1};

    \draw[green, thick] (1) -- (2);
    \draw[green, thick] (2) -- (3); 
    \draw[black, thick] (3) -- (4);
    \draw[blue, thick] (4) -- (5);

    \draw[green, thick] (6) -- (7);
    \draw[green, thick] (7) -- (8); 
    \draw[black, thick] (8) -- (9);
    \draw[blue, thick] (9) -- (10);

    \draw[->, blue, thick, bend left=45] (5) to (1);
    \draw[->, blue, thick, bend left=45] (10) to (6);

    \draw[->, thick] (-1,0) -- (1,0);
    \node[minimum size=0pt,inner sep=0pt,label=below: Figure 12] (name) at (0,-0.7){};
    \end{tikzpicture}
\end{eqnarray*}

We have now examined all possible configurations involving 5 pebbles on $C_5$ and demonstrated that each is $v$-solvable. With this, our proof is complete.

\end{proof}

Hence, we obtain:

\begin{Corollary}
For a cycle $C_5$, the pebbling number is given by:

\[\pi(C_5) = 5.\]
\end{Corollary}

To continue completing the proof for the Petersen graph, we need to establish both inequalities: $\pi(P) > 9$ and $\pi(P) \leq 10$. For the first inequality, $\pi(P) > 9$, we will demonstrate that there exists a configuration with 9 pebbles that is $r$-unsolvable for some vertex $r$. Similar to the example used for $C_5$, we place one pebble on each of the 9 vertices of the Petersen graph, leaving one vertex without a pebble. This configuration will show that the graph is not solvable with 9 pebbles, thereby proving that $\pi(P) > 9$.

Then we will move on to prove $\pi(P)\leq 10$. According to Figure 9, we say that we have an outer cycle and an inner cycle, while the entire graph is symmetric and can be rotated. It is straightforward to divide all possibilities into two categories: the root being on the outer cycle or the inner cycle. However, since we know that the outer cycle and inner cycle are completely equivalent, it is sufficient for us to analyze just one of these cycles to cover all possibilities. Thus, we say the root $r$ is on the outer cycle. Due to the Petersen graph's property, every vertex has exactly three neighbours. Now, we have three cases:

\begin{itemize}
    \item \textbf{Case 1:} One of $r$'s neighbours has two pebbles which directly proves that it is $r$-solvable.
  
    \item \textbf{Case 2:} One of $r$'s neighbours has one pebble. We denote this neighbour as $s$. Consequently, we have two new separate cycles with 5 vertices each, including $r$ or $s$ (See Figure 13). We denote the two sub-cycles formed from splitting the Petersen graph as $C_5^r$ and $C_5^s$, respectively. 

       \begin{eqnarray*}
    \begin{tikzpicture}[scale=0.7]
    \tikzstyle{every node}=[thick,minimum size=2pt, inner sep=1pt]
    
    \node[circle, fill=red] (1) at (-7,2){};
    \node[circle, fill=blue] (2) at (-8.9,0.6){};
    \node[circle, fill=black] (3) at (-8.2,-1.6){};
    \node[circle, fill=black] (4) at (-5.8,-1.6){};
    \node[circle, fill=black] (5) at (-5.1,0.6){};

    \node[circle, fill=black] (6) at (-7,0.8){};
    \node[circle, fill=black] (7) at (-7.6,0.5){};
    \node[circle, fill=black] (8) at (-7.6,-0.5){};
    \node[circle, fill=black] (9) at (-6.4,-0.5){};
    \node[circle, fill=black] (10) at (-6.4,0.5){};

    \node[circle, fill=red] (11) at (0,2){};
    \node[circle, fill=black] (12) at (-1.9,0.6){};
    \node[circle, fill=black] (13) at (-1.2,-1.6){};
    \node[circle, fill=black] (14) at (1.2,-1.6){};
    \node[circle, fill=black] (15) at (1.9,0.6){};

    \node[circle, fill=blue] (16) at (0,0.8){};
    \node[circle, fill=black] (17) at (-0.6,0.5){};
    \node[circle, fill=black] (18) at (-0.6,-0.5){};
    \node[circle, fill=black] (19) at (0.6,-0.5){};
    \node[circle, fill=black] (20) at (0.6,0.5){};

    \node[circle, fill=red] (21) at (7,2){};
    \node[circle, fill=black] (22) at (5.1,0.6){};
    \node[circle, fill=black] (23) at (5.8,-1.6){};
    \node[circle, fill=black] (24) at (8.2,-1.6){};
    \node[circle, fill=blue] (25) at (8.9,0.6){};

    \node[circle, fill=black] (26) at (7,0.8){};
    \node[circle, fill=black] (27) at (6.4,0.5){};
    \node[circle, fill=black] (28) at (6.4,-0.5){};
    \node[circle, fill=black] (29) at (7.6,-0.5){};
    \node[circle, fill=black] (30) at (7.6,0.5){};
    
    \draw [thick](1) -- (2);
    \draw [thick][pink](2) -- (3);
    \draw [thick](3) -- (4);
    \draw [thick][green](4) -- (5);
    \draw [thick][green](5) -- (1);
    \draw [thick](6) -- (8);
    \draw [thick](7) -- (9);
    \draw [thick][pink](8) -- (10);
    \draw [thick][green](9) -- (6);
    \draw [thick][pink](10) -- (7);
    \draw [thick][green](1) -- (6);
    \draw [thick][pink](2) -- (7);
    \draw [thick][pink](3) -- (8);
    \draw [thick][green](4) -- (9);
    \draw [thick](5) -- (10);

    \draw [thick][green](11) -- (12);
    \draw [thick][green](12) -- (13);
    \draw [thick][green](13) -- (14);
    \draw [thick][green](14) -- (15);
    \draw [thick][green](15) -- (11);
    \draw [thick][pink](16) -- (18);
    \draw [thick][pink](17) -- (19);
    \draw [thick][pink](18) -- (20);
    \draw [thick][pink](19) -- (16);
    \draw [thick][pink](20) -- (17);
    \draw [thick](11) -- (16);
    \draw [thick](12) -- (17);
    \draw [thick](13) -- (18);
    \draw [thick](14) -- (19);
    \draw [thick](15) -- (20);

    \draw [thick][green](21) -- (22);
    \draw [thick][green](22) -- (23);
    \draw [thick](23) -- (24);
    \draw [thick][pink](24) -- (25);
    \draw [thick](25) -- (21);
    \draw [thick][green](26) -- (28);
    \draw [thick][pink](27) -- (29);
    \draw [thick](28) -- (30);
    \draw [thick](29) -- (26);
    \draw [thick][pink](30) -- (27);
    \draw [thick][green](21) -- (26);
    \draw [thick](22) -- (27);
    \draw [thick][green](23) -- (28);
    \draw [thick][pink](24) -- (29);
    \draw [thick][pink](25) -- (30);
    
    \node[minimum size=0pt,inner sep=0pt,label=below: Figure 13: Different Selections of $C_5^r$ and $C_5^s$ Cycles in the Petersen Graph.] (name) at (0,-2.5){};
    \end{tikzpicture}   
    \end{eqnarray*}

    Since we have already proven that $\pi(C_5) = 5$, when there are more than 5 pebbles on $C_5^r$, the configuration should be $r$-solvable. On the other hand, we have:

    \[C(C_5^r) \leq 4 \quad \text{and} \quad C(C_5^s) \geq 6.\]

    This implies that, on $C_5^s$, even without moving the pebble from $s$, we can still place another pebble on $s$, ensuring that there will ultimately be 2 pebbles on $s$. Thus, we conclude that with 10 pebbles on the Petersen graph, any configuration is $r$-solvable.

    \item \textbf{Case 3:} There is no pebble on $r$'s neighbours. Now, this suggests that we have 10 pebbles on the other 6 vertices (See Figure 14).
    \begin{eqnarray*}
    \begin{tikzpicture}[scale=0.7]
    \tikzstyle{every node}=[thick,minimum size=2pt, inner sep=1pt]
    
    \node[circle, fill=black] (1) at (0,2){};
    \node[circle, fill=black] (2) at (-1.9,0.6){};
    \node[circle, fill=red] (3) at (-1.2,-1.6){};
    \node[circle, fill=red] (4) at (1.2,-1.6){};
    \node[circle, fill=black] (5) at (1.9,0.6){};

    \node[circle, fill=black] (6) at (0,0.8){};
    \node[circle, fill=red] (7) at (-0.6,0.5){};
    \node[circle, fill=red] (8) at (-0.6,-0.5){};
    \node[circle, fill=red] (9) at (0.6,-0.5){};
    \node[circle, fill=red] (10) at (0.6,0.5){};
    \draw (1) -- (2);
    \draw (2) -- (3);
    \draw (3) -- (4);
    \draw (4) -- (5);
    \draw (5) -- (1);
    \draw (6) -- (8);
    \draw (7) -- (9);
    \draw (8) -- (10);
    \draw (9) -- (6);
    \draw (10) -- (7);
    \draw (1) -- (6);
    \draw (2) -- (7);
    \draw (3) -- (8);
    \draw (4) -- (9);
    \draw (5) -- (10);
    
    \node[minimum size=0pt,inner sep=0pt,label=below: Figure 14: The Petersen Graph with 10 Pebbles Distributed Among Red Vertices.] (name) at (0,-2.5){};
    \end{tikzpicture}   
    \end{eqnarray*}
   In Figure 14, we observe that every neighbour of $r$ is connected to two of the six vertices. This allows us to partition these six vertices into three groups, where the vertices in each group are connected to the same neighbour of $r$. Since there are 10 pebbles in total, at least one of these groups must contain at least 4 pebbles. If both vertices in this group have 2 pebbles each, the configuration is immediately $r$-solvable because we can move two pebbles to the neighbour of $r$ and then one pebble to $r$ itself. In the alternative case, if one vertex in the group has 1 pebble (denote this vertex as $m$) and the other vertex has 3 pebbles (denote this vertex as $n$), the configuration is still solvable. Since $m$ is connected to two other vertices from the other groups, if either of these two vertices has 2 pebbles, we can move the required pebble to $r$. However, if neither of those vertices has 2 pebbles, then the remaining two vertices (outside this group) must contain at least 4 pebbles. In fact, both of these two vertices are connected to $n$, allowing us to move at least one pebble to $n$ and subsequently move one pebble to $r$, completing the configuration as $r$-solvable.
\end{itemize}
Thus we complete our proof.

\end{proof}

So far, we have presented several effective examples of pebbling numbers on specific graphs through straightforward proofs. Over nearly 50 years of research and development, many significant discoveries related to pebbling numbers have been made. For a detailed summary of these contributions, refer to \href{https://www.people.vcu.edu/~ghurlbert/pebbling/pebb.html}{Hurlbert's Pebbling website}. A condensed overview is provided in Table 15 below\cite{Chung,Hurlbert99}.

\renewcommand{\arraystretch}{1.5}
\begin{center}
\begin{tabularx}{1.1\textwidth}{ 
  | >{\raggedright\arraybackslash}p{0.3\textwidth}
  | >{\centering\arraybackslash}X 
  | >{\raggedleft\arraybackslash}X | }

 \hline
 \textbf{The Graph Family} & \textbf{Its Pebbling Number Formula}\\ 
 \hline
 Complete graphs & $\pi(K_n)=n$ \\
 \hline
 Paths  & $\pi(n)=2^{n-1}$  \\
 \hline
 The Petersen Graphs & $\pi(P)=10$ \\
 \hline
 Cycles  & For $k\geq 1$, $\pi(C_{2k})=2^k$ and $\pi(C_{2k+1})=2\left\lfloor (\frac{2k+1}{3}) \right\rfloor+1$  \\
 \hline
 Cubes & $\pi(Q_d)=2^d$\\
 \hline
 Trees  &$\pi(T)=\left( \sum_{i=1}^{m}2^{q_i} \right)-m+1$  \\
 \hline
\end{tabularx}
\end{center}

\begin{center}
\text{Table 15: Pebbling Numbers of Various Graph Families and Their Formulas.}
\end{center}

\section{Hurlbert's Linear Optimization Technique}\label{Sec:Linear Optimization Technique}

Although we have established general results for the pebbling numbers of several simple graphs, computing the pebbling number for an arbitrary graph remains a difficult challenge. Watson \cite{Watson} demonstrated that determining whether a given configuration is solvable is a NP-complete problem. Additionally, Clark and Milans \cite{CM} showed that deciding whether $\pi(G) \leq k$ is a $\Pi_{2}^{p}$-complete problem.

Given this complexity, generating bounds for the pebbling number of a graph can significantly assist in the research of pebbling numbers. Recognizing the complexity of the task, Hurlbert introduced a novel approach using the linear optimization technique to estimate these bounds more efficiently. These Linear Optimization methods have proven effective in refining the bounds for pebble numbers, particularly in large graphs, and have become a state-of-the-art computational strategy in the field. More specifically, Hurlbert presented a linear optimization framework to establish upper bounds on pebbling numbers, providing a cornerstone for future research. The key idea in this optimization problem is the weight function, which we will now explore in detail.

\subsection{The Weight Function Lemma}

For different types of graphs, we need to construct tailored weight functions that are suited to each graph's properties. To ensure a clear and thorough understanding, we will begin by examining a proved Lemma \ref{Path PN} on a simple yet illustrative structure, the path $P_{n}$. Let $P_{n}$ be the path $v_{1}v_{2}\cdots v_{n}$ on $v$ vertices. Before defining the weight function, it is always crucial to first designate a target vertex.

In Corollary \ref{Path PN}, we established that the pebbling number of $P_n$ is $2^{n-1}$. Now, let's examine this result using a weight function approach. According to Lemma \ref{Root Choice}, we have shown that on a path, selecting either $v_1$ or $v_n$ as the root results in the configuration requiring the maximum number of pebbles to be $r$-unsolvable. Furthermore, since $v_1$ and $v_n$ are equivalent--given that the path can be flipped, making the starting vertex the ending vertex--it is logical to focus on $v_1$. Thus, the closer a vertex is to $v_1$, the fewer pebbles are required to make the configuration $v_1$-solvable. To capture this, vertices closer to $v_1$ should be assigned higher weight values, reflecting their greater influence in requiring fewer pebbles. Consequently, we define a weight function $\omega$ on the vertices of $P_n$ such that $\omega(v_{n-i}) = 2^i$.

We can extend this weight function to an entire configuration by defining $\omega(C) = \sum_{v \in V(G)} \omega(v) C(v)$, where $C(v)$ represents the number of pebbles on vertex $v$. The purpose of defining this weight function is to quantify the overall "pebble weight" of a configuration. This function will provide insight into how to determine if a graph is $v_1$-solvable based on its weights. To explain this statement, we will now introduce the following result, building on the weight function defined for a path:

\begin{Lemma}\label{Path total weight}
For any $v_1$-unsolvable configuration of pebbles on $V(P_n)$, the total weight $$\omega(C)=\sum_{i=2}^{n} \omega(v_i) C(v_i)$$ of the configuration $C$ is at most $2^{n-1}-1$.
\end{Lemma}

\begin{proof}
We prove the lemma by induction on $n$, the length of the path $P_n$.

For $k=2$, placing at most one pebble on $v_2$ can make the configuration  $\mathrm{P}_2$ $ v_1$-unsolvable. In this case, $\omega(\mathrm{P}_2)=1$. 

The induction assumes that for a path length $k\geq2$, the maximum total weight for the $v_1$-unsolvable configuration equals $2^{k-1}-1$. Let $\mathrm{P}_{k+1}$ be a  $v_1$-unsolvable configuration on ${P}_{k+1}$.  Let   $\mathrm{P}_{k}$ be the restriction  of $\mathrm{P}_{k+1}$ to the sub-path $v_3\cdots v_{k+1}$ and $\mathrm{P}_{k}(v_2)=0$. The overall structure is shown below (See Figure 16).

\begin{eqnarray*}
\begin{tikzpicture}[scale=0.8, thick]
   \tikzstyle{every node}=[thick,minimum size=1pt, inner sep=2pt]
    \node[circle, fill=black, minimum size=1pt, label=below:{\footnotesize $v_1$}, label=left:{\footnotesize $\omega(v_1) = 2^k$}] (1) at (0, 0) {};
    \node[circle, fill=black, minimum size=1pt, label=below:{\footnotesize $v_2$}, label=above:{\footnotesize $\omega(v_2) = 2^{k-1}$}] (2) at (2, 0) {};
    \node[circle, fill=black, minimum size=1pt, label=below:{\footnotesize $v_3$}] (3) at (4, 0) {};
    
    \node[label=below:{$\cdots$}] (dots) at (5.5, 0) {};
    
    \node[circle, fill=black, minimum size=1pt, label=below:{\footnotesize $v_k$}, label=above:{\footnotesize $\omega(v_k) = 2$}] (k) at (7, 0) {};
    \node[circle, fill=black, minimum size=1pt, label=below:{\footnotesize $v_{k+1}$}, label=right:{\footnotesize $\omega(v_{k+1}) = 1$}] (k1) at (9, 0) {};

    \draw (1) -- (2);
    \draw (2) -- (3);
    \draw (k) -- (k1);

    \draw[decoration={brace,mirror,raise=15pt},decorate] (2) -- (k1) 
    node[midway, below=20pt] {\footnotesize Path $P_k$};
    
    \node[minimum size=0pt, inner sep=0pt, label=below:Figure 16: Representation of the Path Graph $P_k$ with Weighted Vertices.] (caption) at (4.5, -1.5) {};
\end{tikzpicture}
\end{eqnarray*}

\begin{itemize}
    \item \textbf{Case 1:} If $\mathrm{P}_{k+1}(v_2)=1$. Then,   $\mathrm{P}_{k}$ will be $v_2$-unsolvable. By induction, $\omega(\mathrm{P}_{k})\leq  2^{k-1}-1$ that is $\sum_{i=3}^{k+1} \omega(v_i)\mathrm{P}_{k}(v_i)\leq 2^{k-1}-1$. Since $\mathrm{P}_{k+1}(v_1)=0$, $\mathrm{P}_{k+1}(v_2)=1$ and $\sum_{i=3}^{k+1} \omega(v_i)\mathrm{P}_{k}(v_i)\leq 2^{k-1}-1$, we have:
    \begin{eqnarray*}
       \omega(\mathrm{P}_{k+1}) &=& \sum_{i=2}^{k+1} \omega(v_i) \mathrm{P}_{k+1}(v_i)\\
       &=& \omega(v_2) \mathrm{P}_{k+1}(v_2)+\sum_{i=3}^{k+1} \omega(v_i) \mathrm{P}_{k+1}(v_i)\\
       &=&2^{k-1}+\sum_{i=3}^{k+1} \omega(v_i) \mathrm{P}_{k}(v_i)\\
        &\leq&2^{k-1}+2^{k-1}-1 =2^{k}-1.
    \end{eqnarray*}

    \item \textbf{Case 2:} If $\mathrm{P}_{k+1}(v_2)=0$.Notice that taking pebbling moves will either decrease or preserve the total weight of the configuration. So, the total weight $\omega(\mathrm{P}_{k+1})$ in \textbf{Case 2} will not be more than in \textbf{Case 1}.
\end{itemize}
 Hence, we have $\omega(\mathrm{P}_{k+1})\leq  2^{k}-1$. Clearly,  we can choose $\mathrm{P}_{k+1}(v_1)=0$, $\mathrm{P}_{k+1}(v_i)=1$ for all $i=2,\ldots,k+1$ to make it a maximum total weight  $v_1$-unsolvable configuration. We are done with our proof.

\end{proof}

Now, based on the results above, we can observe that the expression $2^{n-1}-1=\sum_{i=2}^{n}2^{n-i}$ precisely represents the total weight of a configuration in which exactly one pebble is placed on each vertex except the root. We denote this particular configuration as $\mathbf{P}_n$. Therefore, we have:
 
\begin{Corollary}\label{Base case P}
For any $v_1$-unsolvable configuration of pebbles on $V(\mathbf{P}_n)$, the weight function satisfies:

\[\omega(C)\leq \omega(\mathbf{P}_n),\] where $\mathbf{P}_n$ is the configuration on the path $P_n$ that places one pebble on every vertex except the root $v_1$, which has no pebbles.

\end{Corollary}

The formula we derived reveals the purpose of applying the weight function to the configuration: it allows us to transform a complex graph-based game into straightforward numerical calculations. By using an appropriate weight function, we can simplify the analysis of the game, making it more manageable and easier to solve.\\

Based on our understanding of the weight functions as applied to paths, we will now generalize this method and its results into trees, which encompass multiple paths. Let $G$ be a graph and $T$ a subtree(containing at least 2 vertices) of $G$ rooted at vertex $r$. For simplicity, within the subtree $T$, we denote by $v^+$ the parent of a vertex $v\in V(T)$, i.e., the neighbour of $v$ in $T$ that is one step closer to $r$. Correspondingly, we can refer to $v$ as a child of $v^+$ (See Figure 17). \\

\begin{eqnarray*}
    \begin{tikzpicture}[scale=0.7]
    \tikzstyle{every node}=[thick,minimum size=4pt, inner sep=1pt]
    \node[circle, fill=black,label=right:{$r$}] (1) at (0,0){};
    \node[circle, fill=black,label=right:{$v^+$}] (2) at (0,-1.5){};
    \node[circle, fill=black,label=right:{$v$}] (3) at (1,-3){};
    \node[circle, fill=black] (4) at (-1,-3){};
    \node[circle, fill=black] (5) at (-1.5,-1.5){};
    \node[circle, fill=black] (6) at (1.5,-1.5){};
    \draw (1)--(2);
    \draw (3)--(2);
    \draw (4)--(2);
    \draw (5)--(1);
    \draw (6)--(1);
 \node[minimum size=0pt,inner sep=0pt,label=below: Figure 17: Illustration of a Subtree $T$ with Vertex $v$ and Its Parent $v^+$.] (name) at (0,-3.8){};
    \end{tikzpicture}   
    \end{eqnarray*}
    
Next, we will apply a weight function $\omega$ to the subtree $T$, which we refer to as a strategy $S(T, r, \omega)$.
\begin{Definition}
    The strategy $S(T, r, \omega)$ of a rooted tree $(T,r)$ is a triple $(T,r,\omega)$ including a rooted tree $(T,r)$ and a weight function $\omega$ on $(T,r)$ which satisfies a fundamental property: 
    $$\omega(r) = 0, \omega(v^+) = 2\omega(v) \text{ for every } v\in V(G).$$
\end{Definition}

  Similar to the configuration $\mathbf{P}$ on paths, we now canonical define a configuration $\mathbf{T}$ where $C(r) = 0$, $C(v) = 1$ for all $v \in V(\mathbf{T})$, and $C(v) = 0$ elsewhere, representing the number of pebbles on each vertex.

Keeping the notations above, we will now give the most important lemma of weight functions \cite{Hurlbert}:

\begin{Lemma}[Weight Function Lemma]\label{WFL}

Let $S(T, r, \omega)$ be a strategy applied on a subtree $T$ of $G$ rooted at $r$, with associated weight function $\omega$. Suppose that $C$ is an $r$-unsolvable configuration of pebbles on $V(G)$. Then $\omega(C)\leq \omega({\mathbf{T}})$.

\end{Lemma}

\begin{proof}

We use proof of contradiction and induction here. The base case is when $T$ is a path, which we already proved as Corollary \ref{Base case P}.\\

For other cases, assume $\omega(C) > \omega(\mathbf{T})$. Start by selecting a leaf of $T$, denoted by $y$. Define $P$ as the path from the leaf $y$ to the root $r$ in $T$. Let $P_y$ represent the subpath from $y$ to the nearest vertex $x$ on $P$ that has a degree of at least 3 in $T$ (or $r$ if no such vertex exists). In other words, $x$ is the last vertex on $P$ that is connected to another vertex outside of $P$, or $x$ is the root $r$ if no such vertex exists. Now, define $T'$ as the tree obtained by removing $P_y$ from $T$ and reattaching the vertex $x$. This tree $T'$ can be regarded as an equivalent subtree of $T$ (See Figure 18).

\begin{eqnarray*}
    \begin{tikzpicture}[scale=0.8]
    \tikzstyle{every node}=[thick,minimum size=4pt, inner sep=1pt]
    
    \node[circle, fill=black, label=above:$r$] (r) at (0,6) {};
    \node[circle, fill=black, label=right:$x$] (x) at (0,4.5) {};
    \node[circle, fill=black, label=left:$y$] (y) at (-3.5,1.5) {};
    \node[circle, fill=black] (y1) at (-2,3) {};
    \node[circle, fill=black] (y2) at (2,3) {};
    \node[circle, fill=black] (z2) at (0.5,1.5) {};
    \node[circle, fill=black] (z3) at (2,1.5) {};
    \node[circle, fill=black] (z4) at (2,0) {};
    
    \draw (r)--(x);
    \draw (x)--(y1);
    \draw (x)--(y2);
    \draw (y1)--(y);
    \draw (y2)--(z2);
    \draw (y2)--(z3);
    \draw (z3)--(z4);
    
    \draw[red,thick] (y)--(y1)--(x);
    
    \draw[blue,thick] (x)--(y2);
    \draw[blue,thick] (y2)--(z2);
    \draw[blue,thick] (y2)--(z3);
    \draw[blue,thick] (z3)--(z4);

    \node[below, align=center] at (0,-1) {Figure 18: Tree $T$ with path $P_y$ highlighted in red, and subtree $T'$ in blue.};
    
    \end{tikzpicture}
\end{eqnarray*}

Among all $r$-unsolvable configurations, let $C$ be the configuration with the largest weight on $T'$. The restriction of the weight function $\omega$ to $T'$, denoted $\omega'$, still represents a valid strategy for the root $r$ we have not changed the strategy, only reduced the size of the tree. Therefore, we have the inequality:

\[
\omega'(C) \leq \omega'(T').
\]

Similarly, by restricting the weight function $\omega$ to the path $P_y$, denoted by $\omega_y$, we confirm that $P_y$ remains a valid strategy for the root $x$. Thus, we have:
\[
\omega(C) = \sum_{\substack{v \in T \\ v \neq r}} \omega(v) C(v).
\]
Given that $T' + P_y = T - r$ according to the structure of the graph, and knowing that $\omega(r) = 0$, we can express:
\[
\omega(C) = \sum_{\substack{v \in T \\ v \neq r}} \omega(v) C(v) = \sum_{\substack{v \in T' \\ v \neq r}} \omega(v) C(v) + \sum_{\substack{v \in P_y \\ v \neq r}} \omega(v) C(v).
\]
Since the weight function $\omega$ is consistent across $T$, $T'$, and $P_y$, we obtain:
\[
\sum_{\substack{v \in T' \\ v \neq r}} \omega(v) C(v) + \sum_{\substack{v \in P_y \\ v \neq r}} \omega(v) C(v) = \sum_{\substack{v \in T' \\ v \neq r}} \omega'(v) C(v) + \sum_{\substack{v \in P_y \\ v \neq r}} \omega_y(v) C(v),
\]
which simplifies to:
\[
\omega(C) = \omega'(C) + \omega_y(C).
\]
Similarly, for the entire tree $T$, we have:
\[
\omega(T) = \omega'(T') + \omega_y(P_y).
\]
Recall that, at the very beginning, we assumed:
\[\omega(C) > \omega(T).\]
Expanding both sides, we obtain:
\[
\omega'(C) + \omega_y(C) > \omega'(T') + \omega_y(P_y).
\]
Since we have already established that:
\[
\omega'(C) \leq \omega'(T'),
\]
it follows that:
\[
\omega_y(C) > \omega_y(P_y).
\]
Since we know that $x \neq r$, let $C_x$ be the configuration obtained after moving a pebble from $P_{y} - x$ to $x$. Because $C$ is unsolvable, $C_x$ must also be unsolvable. We then have $\omega(C_x) = \omega(C)$, but at the same time, $\omega'(C_x) > \omega'(C)$, which leads to a contradiction.

Therefore, the original lemma is proven.

\end{proof}

\subsection{The Linear Optimization Technique}

So far, we have provided a fundamental understanding of how weight function and the Pebbling Game relate to each other. Now, we will discuss the Linear Optimization Technique. In general, for the Linear Optimization Technique, we require three basic variables:

\begin{itemize} 
    \item [(a)] A graph $G$, such as trees $T$.
    
    \item [(b)] A configuration of pebbles $C$ on $G$, such as $\mathrm{P}$ or $\mathrm{T}$.

    \item [(c)] A Weight Function $\omega$.
  
\end{itemize}

According to Lemma \ref{WFL}, let $\mathcal{T}$ be the set of all $r$-strategies in a graph $G$ with root vertex $r$ (i.e., $\mathcal{T}$ includes all the subtrees with root $r$ that can be found in $G$). For the optimal value of the integer linear optimization, denoted by $z_{G,r}$, we have:

\[z_{G,r}=\text{Max.}\sum_{v\neq r}C(v) \quad \text{s.t. }\omega(C)\leq \omega(\mathbf{T}) \text{ with a witnessing weight function }\omega.\]

This leads to the following corollary:

\begin{Corollary}

For every graph $G$ and root $r$, we have $\pi(G, r) \leq z_{G,r} + 1$.

\end{Corollary}

\begin{proof}

According to the results of the Linear Optimization Technique, we have:

\[W_{1}=\left\{ {C\mid \text{r-unsolvable}}\right\} \subseteq W_{2}=\left\{ {C\mid \omega(C)\leq \omega(\mathbf{T}})\right\}.\]

This implies:

\[\max_{C\in W_{1}}\sum_{v\neq r}C(v)\leq \max_{C\in W_{2}}\sum_{v\neq r}C(v).\]

Since

\[\pi(G,r)=1+\max_{C\in W_{1}}\sum_{v\neq r}C(v)\]

and

\[z_{G,r}=\max_{C\in W_{2}}\sum_{v\neq r}C(v),\]

we obtain:

\[\pi(G,r)=1+\max_{C\in W_{1}}\sum_{v\neq r}C(v)\leq \max_{C\in W_{2}}\sum_{v\neq r}C(v)+1=z_{G,r}+1. \]

\end{proof}

\subsection{Applications}

Until now, we have provided an overview of the fundamental definitions and properties of the Linear Optimization Technique. Next, we will explore its applications. We will continue to use mathematical methods to demonstrate how the technique operates.

Consider a graph with strategies $S\left(T_1, r, \omega_1\right)$, $S\left(T_2, r, \omega_2\right)$, \dots, $S\left(T_k, r, \omega_k\right)$. For each strategy, according to Lemma \ref{WFL}, if the configuration is $r$-unsolvable, we have:

\[\omega(C) \leq \omega(\mathrm{T}),\] where $C$ represents the configuration and $\mathrm{T}$ represents a unique configuration on a tree, as mentioned in the previous section.

To make the graph $r$-unsolvable, we must ensure that each strategy follows the rule above. By summing both sides, we obtain:

\[\sum_{i=1}^{k} \left( \sum_{v \neq r} \omega_i(v) C(v) \right) \leq \sum_{i=1}^{k} \left( \sum_{v \neq r} \omega(v) \right).\]

Rewriting the inequality, we have:

\[\sum_{i=1}^{k} \left( \sum_{v \neq r} \omega_i(v) C(v) \right) = \sum_{v \neq r} \left( \left( \sum_{i=1}^{k} \omega_i(v) \right) C(v) \right) \leq \sum_{i=1}^{k} \left( \sum_{v \neq r} \omega_i(v) \right).\]

Now we let $\kappa$ represent the minimum sum of weights over all vertices:

\[\kappa = \min\left\{ \sum_{i=1}^{k} \omega_i(v_1), \sum_{i=1}^{k} \omega_i(v_2), \dots, \sum_{i=1}^{k} \omega_i(v_n) \right\} \subseteq \mathbb{Z}.\]

Under this condition, we let $\chi$ represent the total sum of weighted pebbles:

\[\kappa \sum_{v \neq r} C(v) \leq \sum_{v \neq r} \left( \sum_{i=1}^{k} \omega_i(v) \right) C(v) \leq \sum_{i=1}^{k} \left( \sum_{v \neq r} \omega_i(v) \right) = \chi.\]

Thus, we arrive at the following:

\[\sum_{v \neq r} C(v) \leq \frac{\chi}{\kappa}.\]

As our goal is to ensure the configuration is $r$-solvable, we conclude:

\[\pi(G) \leq \frac{\chi}{\kappa} + 1.\]

Therefore, the key to applying the Linear Optimization Technique is determining the value of $\kappa$. However, due to the complexity of graph pebbling, finding the optimal weight function and, in turn, the value of $\kappa$ can be time-consuming. Nonetheless, this technique remains a relatively effective method for generating bounds on the pebbling number.

We now present examples of graphs to which the Linear Optimization Technique can be applied. As a starting point, recall that we have already established the formula for the pebbling number of the Petersen graph, as stated below:

\begin{Theorem}\cite{Hurlbert}

\[\pi(P) = 10,\] where $P$ denotes the Petersen graph. 

\end{Theorem}

From Lemma \ref{pi bound}, we know that $\pi(P) \geq 10$. Next, we will demonstrate how the Linear Optimization Technique can be employed to prove the following lemma:

\begin{Lemma} 

\[\pi(P) \leq 10,\] where $P$ denotes the Petersen graph. 

\end{Lemma}

\begin{proof}

We first establish three distinct strategies for the Petersen graph (See Figure 19).

\begin{eqnarray*}
    \begin{tikzpicture}[scale=0.7] 
    \tikzstyle{every node}=[thick,minimum size=2pt, inner sep=1pt]
    
    \node[circle, fill=red,label=above:{\footnotesize r}] (1) at (-8,2){};
    \node[circle, fill=black,label=left:{\footnotesize 4}] (2) at (-9.9,0.6){};
    \node[circle, fill=black,label=left:{\footnotesize 2}] (3) at (-9.2,-1.6){};
    \node[circle, fill=black,label=right:{\footnotesize 1}] (4) at (-6.8,-1.6){};
    \node[circle, fill=black] (5) at (-6.1,0.6){};

    \node[circle, fill=black] (6) at (-8,0.8){};
    \node[circle, fill=black,label=above:{\footnotesize 2}] (7) at (-8.6,0.5){};
    \node[circle, fill=black,label=right:{\footnotesize 1}] (8) at (-8.6,-0.5){};
    \node[circle, fill=black,label=right:{\footnotesize 1}] (9) at (-7.4,-0.5){};
    \node[circle, fill=black,label=above:{\footnotesize 1}] (10) at (-7.4,0.5){};

    \node[circle, fill=red,label=above:{\footnotesize r}] (11) at (-2,2){};
    \node[circle, fill=black] (12) at (-3.9,0.6){};
    \node[circle, fill=black,label=left:{\footnotesize 1}] (13) at (-3.2,-1.6){};
    \node[circle, fill=black,label=right:{\footnotesize 1}] (14) at (-0.8,-1.6){};
    \node[circle, fill=black] (15) at (-0.1,0.6){};

    \node[circle, fill=black,label=right:{\footnotesize 4}] (16) at (-2,0.8){};
    \node[circle, fill=black,label=above:{\footnotesize 1}] (17) at (-2.6,0.5){};
    \node[circle, fill=black,label=right:{\footnotesize 2}] (18) at (-2.6,-0.5){};
    \node[circle, fill=black,label=right:{\footnotesize 2}] (19) at (-1.4,-0.5){};
    \node[circle, fill=black,label=below:{\footnotesize 1}] (20) at (-1.4,0.5){};

    \node[circle, fill=red,label=above:{\footnotesize r}] (21) at (4,2){};
    \node[circle, fill=black] (22) at (2.1,0.6){};
    \node[circle, fill=black,label=left:{\footnotesize 1}] (23) at (2.8,-1.6){};
    \node[circle, fill=black,label=right:{\footnotesize 2}] (24) at (5.2,-1.6){};
    \node[circle, fill=black,label=right:{\footnotesize 4}] (25) at (5.9,0.6){};

    \node[circle, fill=black] (26) at (4,0.8){};
    \node[circle, fill=black,label=above:{\footnotesize 1}] (27) at (3.4,0.5){};
    \node[circle, fill=black,label=left:{\footnotesize 1}] (28) at (3.4,-0.5){};
    \node[circle, fill=black,label=left:{\footnotesize 1}] (29) at (4.6,-0.5){};
    \node[circle, fill=black,label=above:{\footnotesize 2}] (30) at (4.6,0.5){};
    
    \draw[thick][red] (1) -- (2);
    \draw[thick][red] (2) -- (3);
    \draw[thick][red] (3) -- (4);
    \draw[thick] (4) -- (5);
    \draw[thick] (5) -- (1);
    \draw[thick] (6) -- (8);
    \draw[thick][red] (7) -- (9);
    \draw[thick] (8) -- (10);
    \draw[thick] (9) -- (6);
    \draw[thick][red] (10) -- (7);
    \draw[thick] (1) -- (6);
    \draw[thick][red] (2) -- (7);
    \draw[thick][red] (3) -- (8);
    \draw[thick] (4) -- (9);
    \draw[thick] (5) -- (10);

    \draw[thick] (11) -- (12);
    \draw[thick] (12) -- (13);
    \draw[thick] (13) -- (14);
    \draw[thick] (14) -- (15);
    \draw[thick] (15) -- (11);
    \draw[thick][red] (16) -- (18);
    \draw[thick][red] (17) -- (19);
    \draw[thick][red] (18) -- (20);
    \draw[thick][red] (19) -- (16);
    \draw[thick] (20) -- (17);
    \draw[thick][red] (11) -- (16);
    \draw[thick] (12) -- (17);
    \draw[thick][red] (13) -- (18);
    \draw[thick][red] (14) -- (19);
    \draw[thick] (15) -- (20);

    \draw[thick] (21) -- (22);
    \draw[thick] (22) -- (23);
    \draw[thick][red] (23) -- (24);
    \draw[thick][red] (24) -- (25);
    \draw[thick][red] (25) -- (21);
    \draw[thick] (26) -- (28);
    \draw[thick] (27) -- (29);
    \draw[thick][red] (28) -- (30);
    \draw[thick] (29) -- (26);
    \draw[thick][red] (30) -- (27);
    \draw[thick] (21) -- (26);
    \draw[thick] (22) -- (27);
    \draw[thick] (23) -- (28);
    \draw[thick] (24) -- (29);
    \draw[thick] (25) -- (30);
    
    \node[minimum size=0pt,inner sep=0pt,label=below: Figure 19: Three distinct strategies applied to the Petersen graph.] (name) at (-1.5,-2.5){};
    \end{tikzpicture}   
\end{eqnarray*}

From this choice of strategies, we have:

\[\kappa = \min\left\{ \sum_{i=1}^{k} \omega_i(v_1), \sum_{i=1}^{k} \omega_i(v_2), \dots, \sum_{i=1}^{k} \omega_i(v_n) \right\}=4.\]

Also, we have:

\[\chi=\sum_{i=1}^{k} \left( \sum_{v \neq r} \omega_i(v) \right)=36.\]

Thus, we obtain:

\[\pi(G) \leq \frac{\chi}{\kappa} + 1=10.\]

\end{proof}

By combining these results, we will have a complete proof that the pebbling number of the Petersen graph is exactly 10. Through a clear contrast, we can see that when we find the most appropriate weight function, determining the upper bound of a graph's pebbling number becomes much easier. This process helps us to further accurately confirm the exact pebbling number of the graph. Having established the pebbling number for the Petersen graph, we now extend our analysis to a more complex structure: the fourth weak Bruhat graph.

The 4th weak Bruhat graph $B_4$ is a combinatorial structure associated with the symmetric group $S_4$, which consists of all permutations of four elements (See Figure 20). The graph has 24 vertices, corresponding to the 24 elements of $S_4$, and the edges represent the covering relations in the weak Bruhat order. 

\begin{figure*}[ht]
\centering
\includegraphics[scale=0.5]{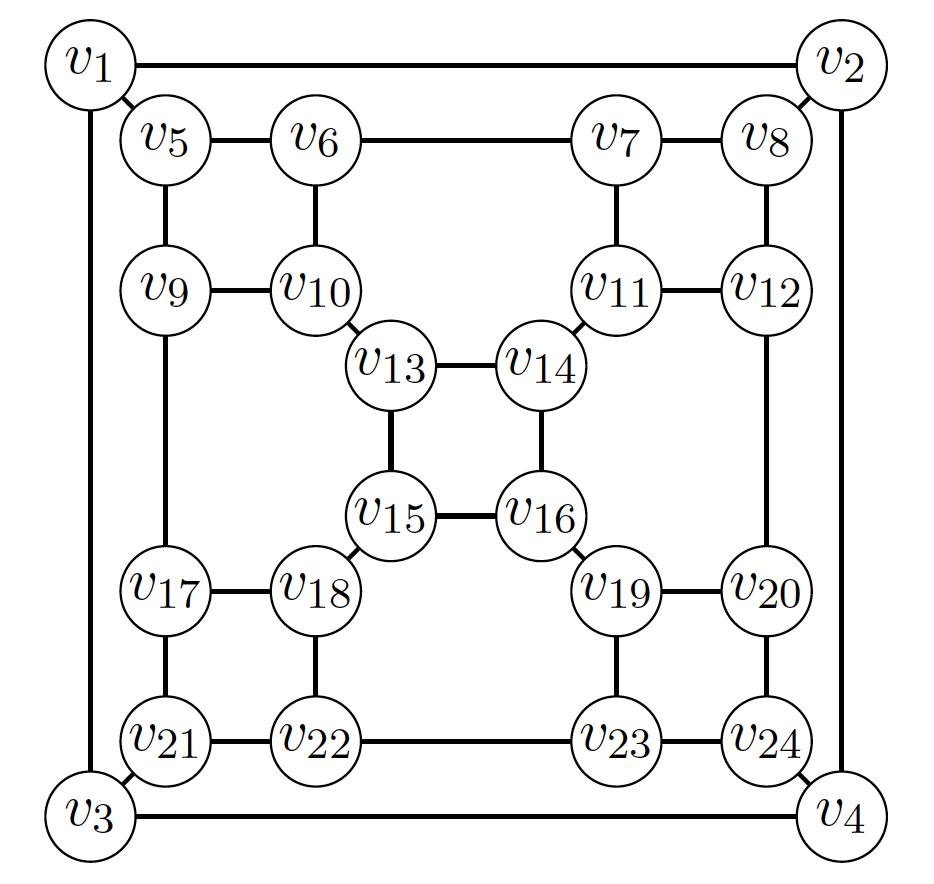}
\text{Figure 20: The Structure of the Fourth Weak Bruhat Graph on 24 Vertices.}
\label{fig:20}
\end{figure*}

Drawing on the construction of such graph structures by Diminic, Jonad, and Carl \cite{AWF}, we can derive the structure of this type of graph (See Figure 20) as well as the application of Linear Optimization Techniques on this graph (See Figure 21). We have:

\begin{Theorem}\label{Bruhat}

Let $B_4$ be the Bruhat graph of order $4$. Then $\pi(B_4)\leq 66.$

\end{Theorem}

\begin{proof}

Similar to the proofs discussed above, we will apply different strategies to the graph $B_4$ (see Figure 21).

Based on the results we have demonstrated, we can now derive the following:

\[\kappa = \min\left\{ \sum_{i=1}^{k} \omega_i(v_1), \sum_{i=1}^{k} \omega_i(v_2), \dots, \sum_{i=1}^{k} \omega_i(v_n) \right\} = 6,\]
and
\[\chi = \sum_{i=1}^{k} \left( \sum_{v \neq r} \omega_i(v) \right) = 395.\]

Thus, we can compute:

\[\pi(G) \leq \frac{\chi}{\kappa} + 1 \approx 66.8.\]

Therefore, we conclude:

\[\pi(G) \leq 66.\]

\end{proof}

\begin{figure*}[ht]
\centering
\includegraphics[scale=0.5]{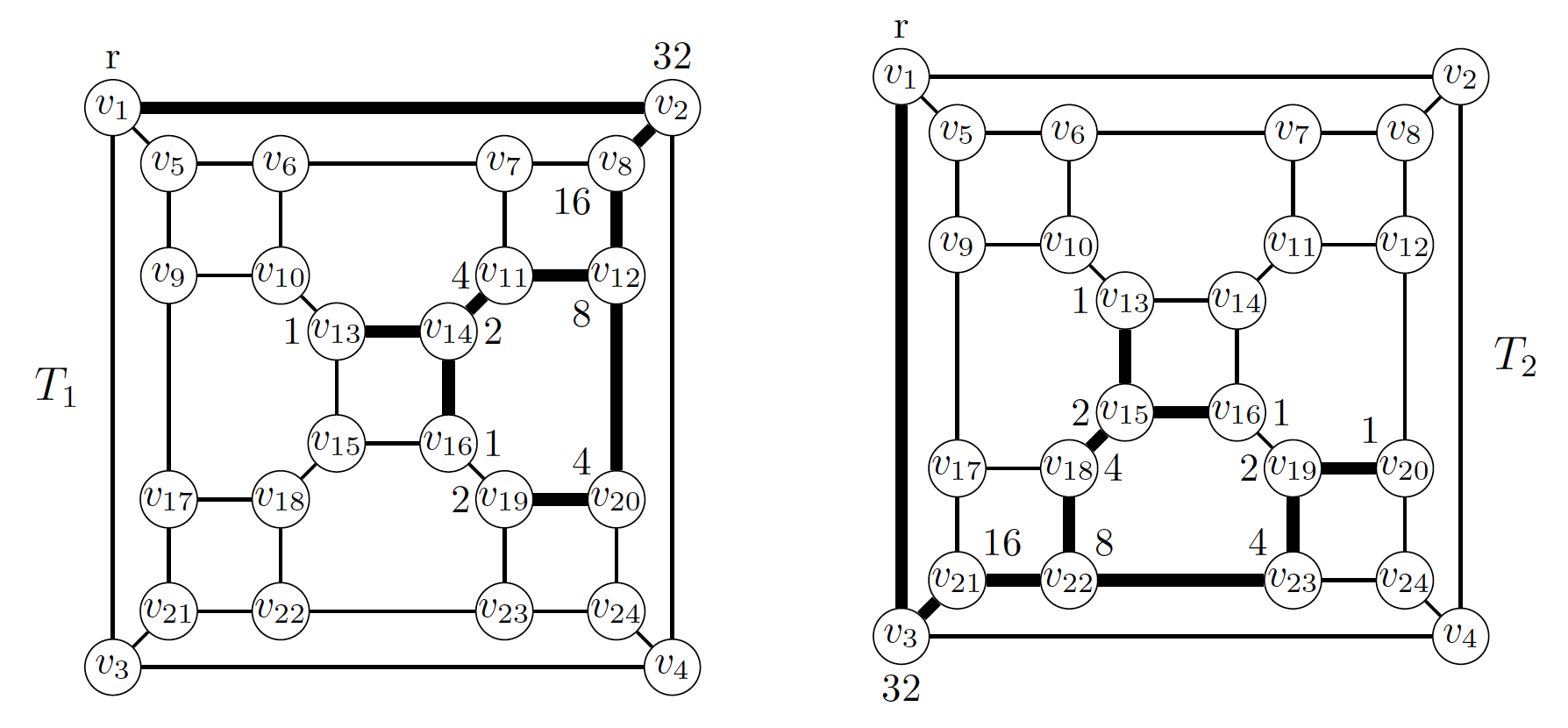}
\centering
\includegraphics[scale=0.5]{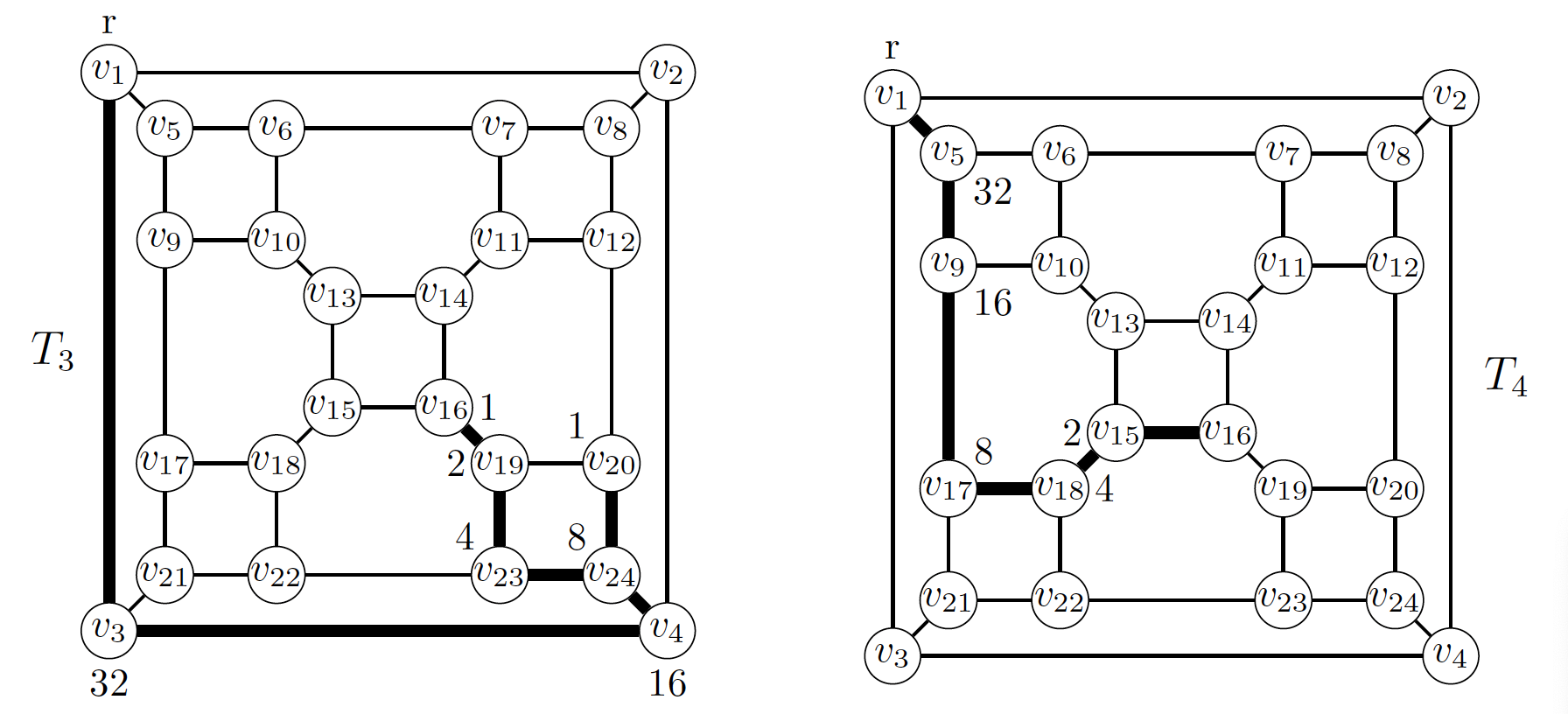}
\centering
\includegraphics[scale=0.5]{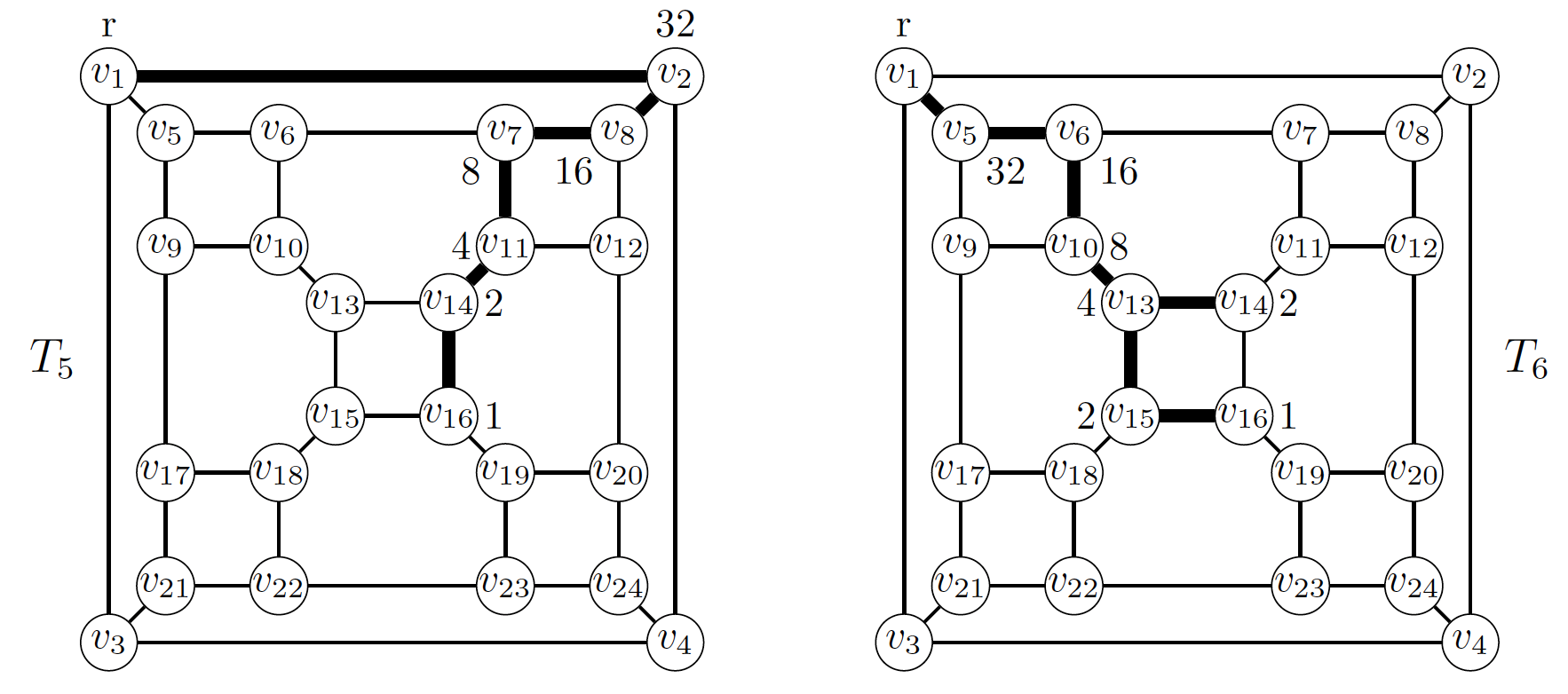}
\text{Figure 21: Set of strategies applied on the Fourth Weak Bruhat Graph with 24 vertices.}
\end{figure*}

After completing the proof for the Petersen graph and applying the Linear Optimization Technique to the Bruhat graph, we now turn our attention to a different class of graphs: trees. Trees, which are acyclic and connected, have unique structural properties that allow us to derive specific results for their pebbling numbers. The following theorem provides a formula for the pebbling number of a tree, rooted at a vertex $r$, referring to results from Chung \cite{Chung} and Hurlbert \cite{Hurlbert}.

\begin{Theorem}
Let $T$ be a tree with a root $r$. The pebbling number $\pi(T, r)$ is given by:
  
  \[\pi(T, r)=\sum_{P \in \mathcal{P}} 2^{e_P}-|\mathcal{P}|+1,\] where $e_P$ denotes the length (number of edges) of the path $P$, and $\mathcal{P}$ represents the set of edge-disjoint paths whose union forms $T$. 
\end{Theorem}

\begin{proof}

Let $C^*$ be a configuration where $2^{e_P} - 1$ pebbles are placed on the leaf $y_P$ of each path $P \in \mathcal{P}$. 

We first show that this configuration $C^*$ is $r$-unsolvable. For each path $P \in \mathcal{P}$, by Lemma~\ref{Path PN}, the configuration is unsolvable for the root of $P$ when restricted to $P$. Since the paths are edge-disjoint, the root $r$ of the entire tree is also unsolvable. Since the total number of pebbles in $C^*$ is $\sum_{P \in \mathcal{P}} 2^{e_P} - |\mathcal{P}|$, we have:

\[\pi(T,r)\geq\sum_{P \in \mathcal{P}} 2^{e_P}-|\mathcal{P}|+1.\]

Next, we show that for any $r$-unsolvable configuration $C$, we have $|C| \leq \sum_{P \in \mathcal{P}} 2^{e_P} - |\mathcal{P}|$, thus proving that $C^*$ is optimal by induction.

For $|\mathcal{P}| = 1$, $T$ is a path, and by Lemma~\ref{Path total weight}, $C^*$ is optimal. Now assume the statement holds for $|\mathcal{P}| = k \geq 2$. For $|\mathcal{P}| = k+1$, let $C^*$ denote the optimal configuration. Define $Q$ as the path in $\mathcal{P}$ whose leaf node has the maximum weight. In case of a tie, choose the path with the shortest length. By the Weight Function Lemma~\ref{WFL}, we have:
\[
\sum_{P \in \mathcal{P} \setminus \{Q\}} \omega(y_P)C(y_P) \leq \sum_{P \in \mathcal{P} \setminus \{Q\}} \omega(P).
\]
By the induction assumption, $C(y_P) \leq 2^{e_P} - 1$ for all $P \neq Q$, with equality holding if and only if $C$ is the maximum configuration. Therefore:
\[
w(y_P) C(y_P) \leq \sum_{P \in \mathcal{P}} w(P) - \sum_{P \in \mathcal{P} \setminus \{Q\}} w(y_P) \left(2^{e_P} - 1\right) = w(Q) = w(y_Q) \left(2^{e_Q} - 1\right).
\]
Thus, we have $C(y_Q) \leq 2^{e_Q} - 1$. Hence, in the optimal configuration $C^*$, pebbles must be placed exclusively on the leaves, completing the proof.


\end{proof}

\section{Conclusion}
In this paper, we have introduced the Pebbling Game, a critical concept in number theory, graph theory, and combinatorics. Over the past half-century, notable mathematicians like F. R. K. Chung and B. Clark, along with many others—both cited and uncited in this paper—have substantially advanced the development and exploration of graph pebbling. Their work has significantly facilitated our understanding of this intriguing topic. This paper focuses on the study of pebbling number bounds, a complex and challenging problem often regarded as hard as a $\Pi^p_2$-complete problem in computational complexity theory. To tackle this challenge, we employed G. Hurlbert's Linear Optimization Technique, a powerful tool that effectively helps in establishing upper bounds for the pebbling number of a graph. By incorporating weight functions within this technique, we provide a systematic approach to calculate these bounds.

For further discussion, while the Linear Optimization Technique greatly enhances the efficiency of determining the pebbling number, it remains complicated to construct an appropriate weight function tailored to a specific graph. The difficulty lies in identifying a weight function that aligns with the unique structure and properties of each graph, which can significantly impact the accuracy of the computed bounds. Therefore, our future efforts will focus on acquiring more pebbling numbers by continuously refining and generating their bounds over time. Additionally, we aim to stay abreast of any new developments or breakthroughs in the field of graph pebbling, with the hope of uncovering novel techniques or results that could further improve the process of calculating pebbling numbers.

\section{Acknowledgement}

I would like to express my deepest gratitude to Professor Carl Yerger for his invaluable guidance and support throughout this research. His expertise and encouragement made this work possible.

\end{document}